\documentclass[11pt, letterpaper]{amsart}
\usepackage[left=1in,right=1in,bottom=0.5in,top=0.8in]{geometry}
\usepackage{amsfonts}
\usepackage{amsmath, amssymb}
\usepackage{graphicx}
\usepackage[font=small,labelfont=bf]{caption}
\usepackage{epstopdf}
\usepackage[pdfpagelabels,hyperindex]{hyperref}
\usepackage{xcolor}
\usepackage{amsthm}
\usepackage{float}
\usepackage{pgfplots}
\usepackage{listings}
\usepackage{longtable}
\usepackage{mathrsfs}

\newtheorem{thm}{Theorem}[section]
\newtheorem{cor}[thm]{Corollary}
\newtheorem{prop}[thm]{Proposition}

\theoremstyle{definition}
\newtheorem{defn}[thm]{Definition}

\theoremstyle{remark}
\newtheorem{rem}[thm]{Remark}

\pgfplotsset{compat=1.9}

\newcommand{\RN}[1]{%
  \textup{\uppercase\expandafter{\romannumeral#1}}%
}
\newcounter{x}

\numberwithin{equation}{section}

\author[Chutian Ma]{Chutian Ma}
\email{cma27@jhu.edu}

\title{A Scattering Result of the Radial Defocusing Cubic Schrodinger Equation in the Hyperbolic Space with Rough Data}
\begin{document}
\begin{abstract}
In this paper, we study the defocusing cubic Schr\"{o}dinger equation  on three dimensional hyperbolic space $\mathbb{H}^3$ with radial initial data in the Sobolev Space $H^s(0<s<1)$. Our main result is that the initial value problem is globally wellposed and scatters for $\frac{15}{16}<s<1$. This is an extension of the high-low method of Bourgain and the work of Staffilani and Yu to the three dimensional hyperbolic space.
\end{abstract}
\maketitle

\section{Introduction}

In this paper we study the cubic dofocusing Schr\"{o}dinger equation on the hyperbolic space $\mathbb{H}^3$ with initial data that is radial and lies in the Sobolev space $H^s$, $0<s<1$:
\begin{equation}\label{eq:01}
    \left\{\begin{aligned}
        &iu_t+\Delta_{\mathbb{H}^3}u=|u|^2u\\
        &u(0,x)=u_0(x)
        \end{aligned}
    \right.
\end{equation}
where $u_0$ is radial and in the Sobolev space $H^s$ with $s<1$.\\
Defocusing NLS with power nonlinearity has been widely studied on Euclidean space. Numerous mass critical and energy critical results are proved on Euclidean space (\cite{dodson2012global},\cite{dodson2016global},\cite{dodson2019global}, \cite{colliander2004global},\cite{visan2007defocusing},\cite{ryckman2007global}). Various results regarding the equation with rough data are also available on Euclidean space, including (\cite{bourgain1998scattering},\cite{colliander2002almost},\cite{colliander2004global}, \cite{colliander2007resonant},\cite{dodson2012globallow},\cite{fang2007global}, \cite{guo2006global}). The equation is of particular interest in the Hyperbolic space, due to the fact that the negative curvature of the manifold results in stronger dispersive estimates and Strichartz estimates. Weighted dispersive and Strichartz estimates for Schr\"{o}dinger equation were obtained in \cite{banica2007nonlinear} and \cite{pierfelice2008weighted}. $L^2$ scattering and $H^1$ subcritical scattering results on $\mathbb{H}^3$ were proved in \cite{banica2008scattering}, \cite{anker2009nonlinear} and later $H^1$ critical scattering result on $\mathbb{H}^3$ was proved in \cite{ionescu2012global}. Also \cite{banica2014global} obtained wellposedness and blow-up results for the focusing NLS. Recently, a rough data scattering result on $\mathbb{H}^2$ was proved in \cite{staffilani2020high} using the frequency high-low method. It was shown that the cubic defocusing NLS with radial rough data scatters if the initial data lies beyond $H^\frac{3}{4}$. This paper aims to extend their result to three dimension.\\ 
 The global existence of (\ref{eq:01}) is equivalent to proving the boundedness of the $H^{\frac{1}{2}}$ norm. In the case when $s\geq 1 $, this follows from the energy conservation and mass conservation. The energy $E(u)$ and mass $M(u)$ are defined by
\begin{equation}\label{eq:02}
    E(u)=\frac{1}{2}\int_{\mathbb{H}^3}|\nabla u|^2 dx + \frac{1}{4}\int_{\mathbb{H}^3}|u|^4dx
\end{equation}
\begin{equation}\label{eq:03}
    M(u)=\int_{\mathbb{H}^3}|u|^2dx
\end{equation}
In the case when $s<1$, the energy is not defined. Thus (1.1) lack a conserved quantity that controls the $H^\frac{1}{2}$ norm. To overcome this difficulty, we split the initial data into two halves and prove an almost conservation law for the low frequency evolution. This paper follows the same idea as in \cite{staffilani2020high}. Our main result is the following theorem:
\begin{thm}
The initial value problem (1.1) is globally wellposed and scatters if $s>\frac{15}{16}$.
\end{thm}
Theorem 1.1 will follow directly from Proposition 4.1 by a continuity argument. 
\pagebreak

\section{Preliminary Results}
In this section we present some of the preliminary results we will need later.

\subsection{Strichartz Estimates}

\begin{defn}[Strichart Spaces]
On $\mathbb{H}^d$, We say that (p,q) is an admissible pair if $p,q\geq2$ and $\frac{2}{p}+\frac{d}{q}=\frac{d}{2}$ and $(p,q,d)\neq (2,+\infty,2)$. Denote p',q' to be the dual of p,q. For $\sigma\in\mathbb{R}$, we define the Strichartz norm $\|\cdot\|_{S^\sigma}$ and its dual $\|\cdot\|_{N^\sigma}$ to be
\begin{equation*}
    \|f\|_{S^\sigma}=\sup\{ \||\nabla|^\sigma f\|_{L^p_tL^q_x}\ |\ {(p,q)\ is\ admissible}\}
\end{equation*}
\begin{equation*}
    \|F\|_{N^\sigma}=\inf\{ \||\nabla|^\sigma F\|_{L^{p'}_tL^{q'}_x}\ |\ {(p,q)\ is\ admissible}\}
\end{equation*}
\end{defn}

\begin{thm}[Strichartz Estimate on $\mathbb{H}^d$ in \cite{ionescu2009semilinear}]
Suppose $\frac{2}{p}+\frac{d}{q}\geq\frac{d}{2}$, where (p,q) $(\Tilde{p},\Tilde{q})$ are admissable, 
Then we have the following estimates

\begin{equation*}
    \|e^{it\Delta}f\|_{L^p_tL^q_x}\leq_{d,p,q} \|f\|_{L^2_x}
\end{equation*}

\begin{equation*}
    \|\int_{s<t}e^{i(t-s)\Delta}F(s)ds\|_{L^p_tL^q_x}\lesssim_{d,p,q,\tilde{p},\tilde{q}} \|F\|_{L^{\tilde{p}'}_tL^{\tilde{q}'}_x}
\end{equation*}

\end{thm}
We will also need the following inhomogeneous Strichartz estimate which extends the possible choices of admissible pairs.
\begin{thm}[Inhomogeneous Strichartz Estimate, \cite{dodson2019defocusing}]\label{th:01}
Suppose $d(\frac{1}{2}-\frac{1}{q})<1$, where $\Tilde{p}, p<\infty$, and  $d(\frac{1}{2}-\frac{1}{q})=\frac{1}{p}+\frac{1}{\Tilde{p}}$.
Then we have the following estimate \begin{equation}
    \|\int_{s<t}e^{i(t-s)\Delta}F(s)ds\|_{L^p_tL^q_x}\lesssim_{d,p,\Tilde{p},q} \|F\|_{L^{\Tilde{p}'}_tL^{q'}_x}
\end{equation}

\end{thm}
\begin{rem}
    Theorem \ref{th:01} is proved in \cite{dodson2019defocusing} under the Euclidean settings. However, the same proof is valid on general manifolds.
\end{rem}

\subsection{Local Smoothing Estimates}\hfill \\

The Schrodinger evolution has the following local smoothing effect due to its dispersive nature.
\begin{thm}[Local Smoothing Estimates in $\mathbb{H}^3$, \cite{lawrie2018local} \cite{kaizuka2014resolvent}]
For any $\delta>0$, We have
\begin{equation}
    \|\langle x \rangle^{-\frac{1}{2}-\delta}|\nabla| e^{it\Delta}f\|_{L^2_{t,x}}\lesssim \||\nabla|^\frac{1}{2}f\|_{L^2_x}
\end{equation}
and
\begin{equation}
    \|\langle x \rangle^{-\frac{1}{2}-\delta}|\nabla|\int_{s<t}e^{i(t-s)\Delta}F(s)ds\|_{L^2_{t,x}}\lesssim \|\langle x \rangle^{\frac{1}{2}+\delta}F\|_{L^2_{t,x}}
\end{equation}
\end{thm}

\begin{cor}
For $0\leq \alpha \leq1$, $\frac{1}{r}=\frac{\alpha}{2}+\frac{3(1-\alpha)}{10}$, we have
\begin{equation}
    \|\langle x \rangle^{-(\frac{1}{2}+\delta)\alpha}|\nabla|e^{it\Delta}f\|_{L^r_{t,x}}\lesssim \||\nabla|^{1-\frac{1}{2}\alpha}f\|_{L^2_x}
\end{equation}
\end{cor}
Proof: Interpolate (2.2) with $\|\nabla e^{it\Delta}f\|_{L^{\frac{10}{3}}_{t,x}}\lesssim \|\nabla f\|_{L^2_x}$ yields the result.

\subsection{Heat-Flow Based Frequency Projection}\hfill \\
As in \cite{lawrie2018local}, heat flow based freqency projection operators are defined by
\begin{equation*}
\begin{aligned}
    &P_{\geq s}f=e^{s\Delta}f\\
    &P_{<s}f=f-P_{\geq s}f\\
    &P_sf=(-s\Delta)e^{s\Delta}f
\end{aligned}
\end{equation*}
These operators serve as a substitute of the Littlewood Paley frequency projections in the Euclidean space. Operators of the form $\chi(-\Delta)$, where $\chi$ is a compactly supported bump function, might seem like a more straightforward choice. But unfortunately, the latter are not $L^p$ bounded in general. $P_{\geq s}$ can be viewed roughly as frequency localization to frequency $\leq s^{-\frac{1}{2}}$, whereas $P_s$ is the frequency localization to frequency $\sim s^{-\frac{1}{2}}$, in the sense that they satisfy the following Bernstein inequalities:

\begin{prop}[Bernstein Inequality]
\begin{equation*}
\begin{aligned}
    &\|P_{<s}f\|_{L^p_x}\lesssim s^\frac{1}{2}\|\nabla f\|_{L^p_x}\\
    &\|\nabla P_{\geq s}f\|_{L^p_x} \lesssim s^{-\frac{1}{2}}\|f\|_{L^p_x}\\
    &\|\nabla P_{s}f\|_{L^p_x} \sim s^{-\frac{1}{2}}\|f\|_{L^p_x}
\end{aligned}
\end{equation*}

\end{prop}

\subsection{Radial Sobolev Embeddings on $\mathbb{H}^3$}\hfill \\
We will need the following radial Sobolev inequality.
\begin{prop}[Corollary 2.22 in \cite{staffilani2020high}]
For any $\frac{1}{4}<\alpha<1$ and f radial,
\begin{equation*}
\|sinh(r)f\|_{L^\infty_x(\mathbb{H}^3)}\lesssim \|f\|_{L^2_x(\mathbb{H}^3)}^{1-\frac{1}{4\alpha}}\|(-\Delta)^\alpha f\|_{L^2_x(\mathbb{H}^3)}^{\frac{1}{4\alpha}}
\end{equation*}
\end{prop}
Proof: The proof in \cite{staffilani2020high} can be modified with little effort to hold in $\mathbb{H}^3$, with the radial term $(sinh(r))^\frac{1}{2}$ replaced by $sinh(r)$.

\section{Energy Increment Estimate}

In this section we treat the solution of (1.1) by splitting it into two different evolutions and proving an energy bound for the nonlinear one. Section 3.1 is devoted to local estimates on time interval small enough. Section 3.2 extends the estimate of Section 3.1 to any interval where solution exists. Section 3.3 proves a technical result we needed in Section 3.1.

\subsection{Local Estimates on Small Interval}\hfill \\
In this subsection, we study the solution to (\ref{eq:01}) with initial data $u_0$ on a small interval I, where the scattering norm of the solution is bounded by some small constant
\begin{equation}\label{eq:04}
\|u\|_{L^8_tL^4_x(I)}^8\leq \epsilon
\end{equation}$0<\epsilon=\epsilon(u_0) \ll 1$ is a small constant to be determined later.\\
To begin with, we decompose the initial data $u_0$ by
\begin{equation}
    u_0=P_{<s_0}u_0+P_{\geq s_0}u_0
\end{equation}
where $s_0=s_0(u_0)$ is to be chosen later and let $\psi$ be the linear evolution of $P_{<s}u_0$. In other words, $\psi(t,x)$ is the solution to
\begin{equation} \label{eq:psi}
    \left\{\begin{aligned}
        &i\psi_t+\Delta_{\mathbb{H}^3}\psi=0\\
        &\psi(0,x)=P_{<s_0}u_0
    \end{aligned}
    \right.
\end{equation}
We define $\phi(t,x)$ and v(t,x) to be the solution to the following nonlinear solutions respectively:
\begin{equation} \label{eq:phi}
    \left\{\begin{aligned}
        &i\phi_t+\Delta_{\mathbb{H}^3}\phi=|\phi|^2\phi\\
        &\psi(0,x)=P_{\geq s_0}u_0
    \end{aligned}
    \right.
\end{equation}
and
\begin{equation} \label{eq:v}
    \left\{\begin{aligned}
        &iv_t+\Delta_{\mathbb{H}^3}v=|u|^2u-|\phi|^2\phi\\
        &v(0,x)=0
    \end{aligned}
    \right.
\end{equation}
Write $\zeta=\phi+v$.

Under the assumption that 
\begin{equation}\label{eq:energyassumption}
    E(\phi)\sim s_0^{-(1-s)}
\end{equation}
we have the follwing estimates on the Strichartz norm of $\psi,\phi$ and v.

\begin{prop}
We have the following estimates
\begin{equation}\label{eq:05}
\begin{aligned}
    &\|\psi\|_{S^\sigma(I)}\lesssim s_0^{\frac{1}{2}(s-\sigma)},\ 0\leq\sigma \leq s\\
    &\|\phi\|_{S^\sigma(I)}\lesssim 
        s_0^{-\frac{\sigma}{2}(1-s)},\ 0\leq\sigma\leq 1
\end{aligned}
\end{equation}
In particular, 
\begin{equation}\label{eq:06}
\begin{aligned}
    &\|v\|_{L^8_tL^4_x(I)}\lesssim s_0^{\frac{1}{2}(s-\frac{1}{2})}\\
    &\|\phi\|_{L^8_tL^4_x(I)}\lesssim \epsilon\\
\end{aligned}
\end{equation}
\end{prop}
Proof:
The first inequality in (\ref{eq:05}) follows from Strichartz estimate and Bernstein inequality.\\
Now we prove (\ref{eq:06}) first. Note that
\begin{equation}
    iv_t+\Delta_{\mathbb{H}^3}v=\mathcal{O}(\psi^3+v^3+\phi^2\psi+\phi^2v)
\end{equation}
where we abuse notation since the complex conjugate is not significant in our estimates.
We utilize Thm \ref{th:01} with p=8, q=4
\begin{equation}
\begin{aligned}
    \|v\|_{L^8_tL^4_x}&\lesssim \|\psi\|^3_{L^8_tL^4_x}+\|v\|^3_{L^8_tL^4_x}+\|\phi\|^2_{L^8_tL^4_x}\|\psi\|_{L^8_tL^4_x}+\|\phi\|^2_{L^8_tL^4_x}\|v\|_{L^8_tL^4_x}\\
    &\lesssim s_0^{\frac{3}{2}(s-\frac{1}{2})}+\|v\|^3_{L^8_tL^4_x}+(\epsilon+\|v\|^2_{L^8_tL^4_x})^2s_0^{\frac{1}{2}(s-\frac{1}{2})}+(\epsilon+\|v\|_{L^8_tL^4_x})^2\|v\|_{L^8_tL^4_x}
\end{aligned}
\end{equation}
We get the first inequality in (\ref{eq:06}) after absorbing the terms with small coefficient or higher order exponent into the left hand side. The second inequality in (\ref{eq:06}) follows from the fact that $\phi=u-\psi-v$ and (\ref{eq:04}). With the smallness of $\|\phi\|_{L^8_tL^4_x}$, the second inequality in (\ref{eq:05}) is a straightforward result of Strichartz inequalities.

\begin{cor}[Radial inequalities for $\phi$ and v]
\begin{equation*}
    \begin{aligned}
    &\|sinh(r)\psi\|_{L^\infty_{t,x}}\lesssim s_0^{\frac{1}{2}s-\frac{1}{4}}\\
    &\|sinh(r)\phi\|_{L^\infty_{t,x}}\lesssim s_0^{-\frac{1}{4}(1-s)}
    \end{aligned}
\end{equation*}
\end{cor}
Proof: This is a direct result from proposition 2.8 and the local estimates above.

\begin{prop}[$H^\sigma$ estimate of v when $0\leq\sigma\leq s$]
\begin{equation}
    \|v\|_{S^\sigma(I)}\lesssim s_0^{\frac{1}{2}(\sigma s+s-\sigma-\frac{1}{2})}
\end{equation}
\end{prop}

Proof: Differentiate (3.5) with $|\nabla|^s$, we get
\begin{equation*}
\begin{aligned}
    \||\nabla|^\sigma v\|_{S^0}&\lesssim \||\nabla|^\sigma\psi\|_{S^0}\|\psi\|^2_{L^8_tL^4_x}+\||\nabla|^\sigma v\|_{S^0}\|v\|^2_{L^8_tL^4_x}+
    \||\nabla|^\sigma\phi\|_{S^0}\|\phi\|_{L^8_tL^4_x}\|\psi\|_{L^8_tL^4_x}\\
    &+\||\nabla|^\sigma\psi\|_{S^0}\|\phi\|^2_{L^8_tL^4_x}+\|\phi\|^2_{L^8_tL^4_x}\||\nabla|^sv\|_{S^0}+
    \||\nabla|^\sigma\phi\|_{S^0}\|\phi\|_{L^8_tL^4_x}\|v\|_{L^8_tL^4_x}\\
    &\lesssim s_0^{\frac{1}{2}(\sigma s+s-\sigma-\frac{1}{2})}+o(1)\|v\|_{S^\sigma}
\end{aligned}
\end{equation*}

\begin{prop}[$H^1$ estimate of v]
\begin{equation}
    \|v\|_{S^1(I)}\lesssim s_0^{s-\frac{7}{8}}
\end{equation}
\end{prop}
Let us postpone the proof of prop 3.2 to later. For now, let us assume Proposition 3.4. We are ready to study the solution on longer intervals by cutting the interval into small subintervals, on each of which our local smallness assumption \ref{eq:04} holds.

\subsection{Conditional Global Energy Estimates}

Suppose the solution exists on the interval I and the scattering norm is bounded by some $M>0$
\begin{equation}
\|u\|^8_{L^8_tL^4_x}([0,T])\leq M
\end{equation}
Divide I into subintervals $I_j$ such that on each $I_j$ (\ref{eq:04}) holds. On $I_{j-1}$, evolve $\phi$ and v by the respective equations (\ref{eq:phi}) and (\ref{eq:v}). After the evolution on $I_{j-1}$ is done, we dump the data of v into $\phi$ and evolves $\phi$ with the new initial data on $I_j$, whereas v starts with 0 initial data again. In other words, 
for $t\in (b_j,b_{j+1})$, if we denote $\phi(b_j^+)=lim_{t\downarrow b_j}\phi(t)$ and $\phi(b_j^-)=lim_{t\uparrow b_j}\phi(t)$, $\phi$ and v satisfy
\begin{equation*}
    \left\{\begin{aligned}
        &i\phi_t+\Delta_{\mathbb{H}^3}\phi=|\phi|^2\phi\\
        &\phi(b_j^+,x)=\phi(b_j^-,x)+v(b_j^-,x)
    \end{aligned}
    \right.
\end{equation*}
and
\begin{equation*}
    \left\{\begin{aligned}
        &iv_t+\Delta_{\mathbb{H}^3}v=|u|^2u-|\phi|^2\phi\\
        &v(b_j^+,x)=0
    \end{aligned}
    \right.
\end{equation*}

We can repeat the evolution by the process described above throughout the whole interval I. The local estimates, namely Proposition 3.1-3.4, will hold true on each $I_j$ as long as our assumption (\ref{eq:energyassumption}) remains valid. In other words, the energy increment of $\phi$ resulting from the data of v has to be comparable to $s_0^{-(1-s)}$ in order for our local estimates to remain valid. In fact, under (\ref{eq:energyassumption}),
We have the following conditional global energy increment estimate:
\begin{prop}
Suppose $T>0$ such that $\|u\|^8_{L^8_tL^4_x}([0,T])\leq M$ where $M\sim s_0^{-\frac{1}{2}s+\frac{3}{8}}$. Then we have 
\begin{equation}\label{eq:07}
    E(\phi(T))\leq E(\phi(0))+C\frac{M}{\epsilon}s_0^{\frac{3}{2}s-\frac{11}{8}}\sim s_0^{-(1-s)}
\end{equation}
\end{prop}
Proof:
Divide [0,T] into subintervals $I_j=[b_{j-1},b_j]$ such that $\|u\|^8_{L^8_tL^4_x(I_j)}\leq\epsilon$. The number of such intervals is at most  $\mathcal{O}(\frac{M}{\epsilon})$. During each iteration, the energy increment of $\phi$ is 
\begin{equation}
    \begin{aligned}
    E(\phi(b_j)+v(b_j))-E(\phi(b_{j-1}))&\leq \frac{1}{2}\int|\nabla\phi(b_j)+\nabla v(b_j)|^2dx-\frac{1}{2}\int|\nabla\phi(b_j)|^2dx\\
    &+\frac{1}{4}\int|\phi(b_j)+v(b_j)|^4dx-\frac{1}{4}\int|\phi(b_j)|^4dx\\
    &\leq \int |\nabla\phi(b_j)||\nabla v(b_j)|+|\nabla v(b_j)|^2 dx + \int |\phi(b_j)|^3|v(b_j)|+|\phi(b_j)||v(b_j)^3|dx\\
    &=I+II
    \end{aligned}
\end{equation}

\begin{equation*}
    \begin{aligned}
    I&\leq \|\phi\|_{S^1(I)}\|v\|_{S^1(I)}+\|v\|_{S^1(I)}^2\\
     &\lesssim s_0^{-\frac{1}{2}(1-s)}\cdot s_0^{s-\frac{7}{8}}\\
     &\lesssim s_0^{\frac{3}{2}s-\frac{11}{8}}\\ \\
    II&\leq \|\phi\|_{L^\infty_tL^4_x(I)}^3\|v\|_{L^\infty_tL^4_x(I)}+\|v\|_{L^\infty_tL^4_x(I)}^4\\
    &\lesssim E(\phi)^\frac{3}{4}\|v\|_{L^\infty_tH^\frac{3}{4}_x(I)}+\|v\|_{L^\infty_tH^\frac{3}{4}_x(I)}^4\\
    &\lesssim s_0^{\frac{7}{2}s-\frac{5}{2}} 
    \end{aligned}
\end{equation*}

Combining these terms, 
\begin{equation}
    \Delta E \lesssim s_0^{\frac{3}{2}s-\frac{11}{8}}
\end{equation}

Thus we obtain the local energy increment $\Delta E$ on each interval $I_j$. To justify using the local estimates, we need to ensure that \ref{eq:energyassumption} holds throughout I.
\begin{equation*}
    \frac{M}{\epsilon}\sup_j\Delta E(\phi)\lesssim\frac{M}{\epsilon}s_0^{\frac{3}{2}s-\frac{11}{8}}\lesssim s_0^{-(1-s)}
\end{equation*}
\begin{equation}
    \implies M\sim s_0^{-\frac{1}{2}s+\frac{3}{8}}
\end{equation}
This is the restriction on M arises from.

There are at most $\mathcal{O}(\frac{M}{\epsilon})$ such intervals. Thus \ref{eq:07} holds under the assumption on M.

\subsection{Proof of Proposition 3.4}
It remains to prove Proposition 3.4. 
Take derivative of (3.3),
\begin{equation}
    \left\{  \begin{aligned}
    &i\nabla v_t+ \Delta \nabla v=\mathcal{O}(\nabla\psi\cdot \psi^2+\nabla v\cdot v^2+\nabla\phi\cdot \phi \psi+\nabla\psi \cdot\phi^2+\nabla\phi\cdot \phi v+\nabla v\cdot\phi^2)\\
    &\nabla v(0)=0
    \end{aligned}\right.
\end{equation}

\begin{equation}
    \|v\|_{S^1(I)}\lesssim I+II+III+IV
\end{equation}
where
\begin{equation}
    \begin{aligned}
    &I=\|\nabla\psi\cdot\psi^2\|_{N^0(I)}\\
    &II=\|\nabla v\cdot v^2\|_{N^0(I)}\\
    &III=\|\nabla\phi\cdot\phi\psi\|_{N^0(I)}+\|\nabla \psi\cdot\phi^2\|_{N^0(I)}\\
    &IV=\|\nabla\phi\cdot\phi v\|_{N^0(I)}+\|\nabla v\cdot\phi^2\|_{N^0(I)}
    \end{aligned}
\end{equation}
Estimates for I:
\begin{equation}
    \begin{aligned}
    I\leq \|\langle x \rangle^{-(\frac{1}{2}+\delta)\alpha}\nabla\psi\|_{L^r_{t,x}}
    \|\langle x \rangle^{(\frac{1}{2}+\delta)\alpha}\psi^2\|_{L^{r_1}_{t,x}}\\
    \end{aligned}
\end{equation}
We take $\alpha=\frac{1}{2}$, $\frac{1}{r}=\frac{\alpha}{5}+\frac{3}{10}=\frac{2}{5}$ and $\frac{1}{r_1}=\frac{2}{5}-\frac{\alpha}{5}=\frac{3}{10}$.
Now we apply local smoothing (Corollary 2.4) to estimate the first term in I. Suppose $s\geq \frac{3}{4}$,
\begin{equation}
    \|\langle x \rangle^{-\frac{1}{2}(\frac{1}{2}+\delta)}\nabla\psi\|_{L^r_{t,x}}\lesssim \||\nabla|^{\frac{3}{4}}\psi_0\|_{L^2_x}\lesssim 
    s_0^{\frac{1}{2}(s-\frac{3}{4})}
\end{equation}
For the second term in I, 
\begin{equation}
    \begin{aligned}
     \|\langle x \rangle^{\frac{1}{2}(\frac{1}{2}+\delta)}\psi^2\|_{L^{r_1}_{t,x}(B_1)}
     &\lesssim \|\psi\|^2_{L^{2r_1}_{t,x}}\\
     &\lesssim\|\psi_0\|^2_{H^{\frac{3}{4}}}\\
     &\lesssim s_0^{s-\frac{3}{4}}\\
      \|\langle x \rangle^{\frac{1}{2}(\frac{1}{2}+\delta)}\psi^2\|_{L^{r_1}_{t,x}(B_1^c)}
      &\leq \|sinhr\psi\|_{L^\infty_{t,x}}\|\psi\|_{S^{0}}\\
      &\lesssim s_0^{\frac{1}{2}s-\frac{1}{4}}\cdot s_0^{\frac{1}{2}s}\\
      &\lesssim s_0^{s-\frac{1}{4}}
    \end{aligned}
\end{equation}

Thus \begin{equation}
    I\lesssim s_0^{\frac{1}{2}(s-\frac{3}{4})}\cdot s_0^{s-\frac{3}{4}}\leq s_0^{\frac{3}{2}(s-\frac{3}{4}
    )}
\end{equation}

Estimates for II:
\begin{equation}
    II\leq \|\nabla v\|_{S^1}\|v\|_{L^8_tL^4_x}^2\lesssim o(1)\|v\|_{S^1}
\end{equation}
which is absorbed into the left hand side.

Estimates for III:
We estimate the first term in III in the similar fashion as in I.
Under the asumption that $\frac{1}{2}+\frac{\alpha}{2}\leq s$,
\begin{equation}
\begin{aligned}
        \|\nabla\psi\cdot\phi^2\|_{N^0}&\leq\|\langle x \rangle^{-\frac{1}{2}(\frac{1}{2}+\delta)}\nabla\psi\|_{L^r_{t,x}}\|\langle x \rangle^{\frac{1}{2}(\frac{1}{2}+\delta)}\phi^2\|_{L^{r_1}_{t,x}}\\
        &\lesssim \||\nabla|^{\frac{3}{4}}\psi_0\|_{L^2_x}\left[\|\phi\|^2_{L^{2r_1}_{t,x}(B_1)}+\|sinhr\phi\|_{L^\infty_{t,x}}\|\phi\|_{L^{r_1}_{t,x}(B_1^c)}\right]\\
        &\lesssim \||\nabla|^{1-\frac{1}{2}\alpha}\psi_0\|_{L^2_x}\left[\|\phi\|_{S^1}\phi\|_{L^8_tL^4_x}+\|sinhr\phi\|_{L^\infty_{t,x}}\|\phi\|_{S^0}\right]\\
        &\lesssim s_0^{\frac{1}{2}(s-\frac{3}{4})}[s_0^{-\frac{1}{2}(1-s)}+s_0^{-\frac{1}{4}(1-s)}]\\
        &\lesssim s_0^{s-\frac{7}{8}}
\end{aligned}
\end{equation}
The second term can be estimated by
\begin{equation}
\begin{aligned}
    \|\nabla\phi\cdot\phi\psi\|_{N^0}&\leq\|\nabla\phi\|_{L^\infty_tL^2_x}\|\phi\|_{L^8_tL^4_x}\|\psi\|_{S^{\frac{1}{2}}}\\
    &\lesssim s_0^{-\frac{1}{2}(1-s)}\cdot 1 \cdot s_0^{\frac{1}{2}(s-1)}\\
    &\lesssim s_0^{s-\frac{3}{4}}
\end{aligned}
\end{equation}

Estimates for IV:
The first term in IV
\begin{equation}
    \begin{aligned}
    \|\nabla\phi\cdot\phi v\|_{N^0}&\lesssim \|\nabla\phi\|_{L^2_tL^6_x}\|\phi\|_{L^8_tL^4_x}\|v\|_{L^8_tL^4_x}\\
    &\lesssim s_0^{-\frac{1}{2}(1-s)}\cdot 1 \cdot s_0^{\frac{1}{2}(s-\frac{1}{2})}\\
    &\lesssim s_0^{s-\frac{3}{4}}
    \end{aligned}
\end{equation}
The second term in IV
\begin{equation}
    \begin{aligned}
    \|\nabla v \cdot \phi^2\|_{N^0}&\lesssim \|v\|_{S^1}\|\phi\|_{L^8_tL^4_x}^2\\
    &\lesssim o(1)\|v\|_{S^1}
    \end{aligned}
\end{equation}
which is absorbed into the left hand side.
Combining I to IV, we get $\|v\|_{S^1}\lesssim s_0^{s-\frac{7}{8}}$

\section{Global Wellposedness and Scattering}
We will show that there exists $s_0$ small enough and corresponding $M\sim s_0^{-\frac{1}{2}s+\frac{3}{8}}$, such  that $\|u\|^8_{L^8_tL^4_x[0,T]}\leq M$ for any $[0,T]$ where the solution exists. We will use a bootstrap argument. In fact, we prove the following result

\begin{prop}
For $s>\frac{15}{16}$, there exists $s_0>0$ small enough and $M\sim s_0^{-\frac{1}{2}s+\frac{3}{8}}$, such that the following implication holds:
\begin{equation*}
\|u\|^8_{L^8_tL^4_x[0,T]}\leq M \implies \|u\|^8_{L^8_tL^4_x[0,T]}\leq \frac{1}{2}M
\end{equation*}
\end{prop}

Denote $\zeta=\phi+v$. On [0,T], $\zeta$ satisfies the modified NLS equation
\begin{equation}
    i\zeta_t+\Delta_{\mathbb{H}^3}\zeta=|u|^2u=|\zeta|^2\zeta+\mathcal{N}
\end{equation}
where $\mathcal{N}=|u|^2u-|\zeta|^2\zeta=\mathcal{O}(\psi^3+\psi\zeta^2)$.\\ \\
Due to our choice of M, the result of Proposition 3.3 applies. 
\begin{equation}
    \sup_{t\in[0,T]}E[\zeta(t)]\lesssim s_0^{-(1-s)}
\end{equation}
We estimate the $L^8_tL^4_x$ norm by interpolation
\begin{equation}\label{eq:41}
    \|\zeta\|_{L^8_tL^4_x}^8\leq \|\zeta\|_{L^\infty_tL^4_x}^4\|\zeta\|_{L^4_{t,x}}^4
\end{equation}
$\|\zeta\|_{L^\infty_tL^4_x}$ is bounded by the energy. To estimate $\|u\|_{L^4_{t,x}}$, we utilize the Morawetz inequality developed in \cite{ionescu2009semilinear}

\begin{prop}[Morawetz Estimate on $\mathbb{H}^3$, \cite{ionescu2009semilinear}]
If $\zeta$ satisfies the NLS with error term $\mathcal{N}$,
\begin{equation*}
    i\zeta_t+\Delta_{\mathbb{H}^3}\zeta=|\zeta|^{2}\zeta+\mathcal{N}
\end{equation*}
then the solution $\zeta$ satisfies the following Morawetz estimate
\begin{equation}\label{eq:42}
    \|\zeta\|_{L^{4}_{t,x}}^{4}\lesssim \|\zeta\|_{L^\infty_tL^2_x}\|\zeta\|_{L^\infty_tH^1_x}+\|\mathcal{N}\zeta\|_{L^1_{t,x}}+\|\mathcal{N}\nabla \zeta\|_{L^1_{t,x}}
\end{equation}
\end{prop}

We have the following estimates for the latter two terms in (\ref{eq:42})
\begin{prop}
\begin{enumerate}
    \item[]
    \item \ $\|\mathcal{N}\zeta\|_{L^1_{t,x}}\leq \epsilon\|\zeta\|_{L^4_{t,x}}^4+C_\epsilon s_0^{2s-\frac{1}{2}}$
    \item \ $\|\mathcal{N}\nabla\zeta\|_{L^1_{t,x}}\leq \epsilon\|\zeta\|_{L^4_{t,x}}^4+C_\epsilon s_0^{2s-\frac{3}{2}}$
\end{enumerate}
\end{prop}

Proof:
To see (1), by H\"{o}lder inequality, we have
\begin{equation*}
    \begin{aligned}
    \|\mathcal{N}\zeta\|_{L^1_{t,x}}&\lesssim \|\psi^3\zeta\|_{L^1_{t,x}}+\|\psi\zeta^3\|_{L^1_{t,x}}\\
    &\leq \|\psi\|^3_{L^4_{t,x}}\|\zeta\|_{L^4_{t,x}}+\|\psi\|_{L^4_{t,x}}\|\zeta\|_{L^4_{t,x}}^3\\
    &\lesssim \epsilon\|\zeta\|_{L^4_{t,x}}^4+C_\epsilon\|\psi\|_{L^4_{t,x}}^4
    \end{aligned}
\end{equation*}
Similarly for (2),
\begin{equation*}
    \begin{aligned}
    \|\mathcal{N}\nabla\zeta\|_{L^1_{t,x}}&\lesssim \|\psi^3\nabla\zeta\|_{L^1_{t,x}}+\|\psi\zeta^2\nabla\zeta\|_{L^1_{t,x}}\\
    &\leq \|\psi\|_{L^3_tL^6_x}^3\|\nabla\zeta\|_{L^\infty_tL^2_x}+\|\psi\|_{L^2_tL^\infty_x}\|\zeta\|_{L^4_{t,x}}^2\|\nabla\zeta\|_{L^\infty_tL^2_x}\\
    &\lesssim \epsilon\|\zeta\|_{L^4_{t,x}}^4+C_\epsilon s_0^{2s-\frac{3}{2}}
    \end{aligned}
\end{equation*}
This ends the proof of Proposition 4.2.
\\

Resuming the proof of Proposition 4.1, we plug the results of Proposition 4.3 into (\ref{eq:42}) 
\begin{equation}
\|\zeta\|_{L^4_{t,x}}^4\lesssim s_0^{-\frac{1}{2}(1-s)}+\epsilon\|\zeta\|_{L^4_{t,x}}^4+C_\epsilon s_0^{2s-\frac{1}{2}}+\epsilon\|\zeta\|_{L^4_{t,x}}^4+C_\epsilon s_0^{2s-\frac{3}{2}}
\end{equation}

After absorbing $\epsilon\|\zeta\|_{L^4_{t,x}}^4$ into left hand side, we get 
\begin{equation}
  \|\zeta\|_{L^4_{t,x}}^4\lesssim s_0^{-\frac{1}{2}(1-s)}
\end{equation}
which leads to
\begin{equation}
    \|\zeta\|_{L^8_tL^4_x}^8\lesssim s_0^{-\frac{3}{2}(1-s)}
\end{equation} by interpolation.
And we require the right hand side to be bounded by $\frac{1}{2}M\sim s_0^{-\frac{1}{2}s+\frac{3}{8}}$, i.e.
\begin{equation}
    s_0^{-2(1-s)}\lesssim s_0^{-\frac{1}{2}s+\frac{3}{8}}
\end{equation}
This can be achieved by $s_0$ small enough as long as the exponents satisfy $-\frac{3}{2}(1-s)> -\frac{1}{2}s+\frac{3}{8}$, hence $s>\frac{15}{16}$ is required.

\bibliographystyle{alpha}
\bibliography{NLSH3.bib}


\end{document}